\documentclass{amsart}

\usepackage{amssymb}

\newtheorem{thm}{Theorem}

\newtheorem{cor}[thm]{Corollary}

\theoremstyle{definition}
\newtheorem{defn}[thm]{Definition}

\newtheorem{ques}[thm]{Question}    

\newtheorem{rem}[thm]{Remark}          

\newtheorem{ack}{Acknowledgments}

\newtheorem{defn-thm}[thm]{Definition--Theorem}  
\newtheorem{defn-lem}[thm]{Definition--Lemma}  

\theoremstyle{remark}

\renewcommand{\o}[0]{{\mathcal O}} 

\renewcommand{\a}[0]{{\mathbb A}}

\newcommand{\p}[0]{{\mathbb P}}
\newcommand{\f}[0]{{\mathbb F}}
\newcommand{\q}[0]{{\mathbb Q}}
\newcommand{\map}[0]{\dasharrow}
\newcommand{\qtq}[1]{\quad\mbox{#1}\quad}
\newcommand{\spec}[0]{\operatorname{Spec}}
\newcommand{\pic}[0]{\operatorname{Pic}}

\newcommand{\supp}[0]{\operatorname{Supp}}    
\newcommand{\red}[0]{\operatorname{red}}

\newcommand{\rup}[1]{\lceil{#1}\rceil}

\newcommand{\onto}[0]{\twoheadrightarrow}

\def\into{\DOTSB\lhook\joinrel\to}

\begin{document}
\bibliographystyle{amsalpha}

\title[A conjecture of Ax]{A conjecture of Ax and \\
degenerations of Fano varieties}
\author{J\'anos Koll\'ar}

\today

\maketitle


A field $k$ is called $C_1$ if every homogeneous form
$f(x_0,\dots,x_n)\in k[x_0,\dots,x_n]$
 of degree $\leq n$ has a nontrivial zero.
Examples of $C_1$ fields are finite fields (Chevalley)
and function fields of curves over an algebraically closed  field
(Tsen).

A field is called  {\it  PAC} (pseudo algebraically closed) 
if every geometrically integral $k$-variety has a $k$-point.
An  $k$-variety $X$ is called  geometrically  integral
(or absolutely irreducible) if it is still irreducible and reduced as a 
variety over the algebraic closure $\bar k$.)
These fields were introduced 
in \cite{ax}; see \cite{f-j} for an exhaustive and up to date treatment.

The aim of this paper is to prove 
in characteristic 0 a  conjecture of Ax, posed in \cite[Problem 3]{ax}.

\begin{thm}\label{main.thm.pac} 
 Every PAC field of characteristic 0 is $C_1$.
\end{thm}

\cite[Thm.D]{ax} proves this for fields whose absolute Galois group
is abelian.

Following an idea of \cite{djl}, we deduce 
(\ref{main.thm.pac}) from the next result which
holds for all fields of characteristic zero:

\begin{thm}\label{weakc1.thm}
 Let $k$ be a field of characteristic 0 and
$f_1,\dots,f_s\in k[x_0,\dots,x_n]$ homogeneous polynomials
such that $\sum_i \deg f_i\leq n$.
Let 
$$
X=X(f_1,\dots,f_s):=(f_1=\cdots =f_s=0)\subset \p^n_k
$$
be the subscheme they define in projective $n$-space. Then
\begin{enumerate}
\item  $X$ contains a
geometrically irreducible $k$-subvariety $Y\subset X$.
\item If $k$ is PAC then $X$ has a $k$-point.
\end{enumerate}
\end{thm}

If $k$ is PAC, then $Y$ has a $k$-point which is also a $k$-point of
$X$, thus (\ref{weakc1.thm}.1) implies (\ref{weakc1.thm}.2).
The $s=1$ case of (\ref{weakc1.thm}.2) is precisely (\ref{main.thm.pac}).
The more general version proved here is sometimes called property $C'_1$.

In order to prove (\ref{weakc1.thm}),
we represent
(a subscheme of) the scheme $X(f_1,\dots,f_s)$
as a special fiber of a family $Z\to \p^1$ over the projective line
whose general fiber is a smooth 
hypersurface (resp. complete intersection variety).
The restrictions on the degree are equivalent to
assuming that the canonical class of the general fiber of $Z\to \p^1$
is negative. This approach raises further interesting questions about
degenerations of Fano varieties, we discuss these in
(\ref{src.questions}).

It is thus sufficient  to prove the following
more general result.

\begin{thm}\label{fanofiber.thm}
 Let $k$ be a field of characteristic 0, $C$ a smooth $k$-curve,
$Z$ an irreducible, projective $k$-variety and
$g:Z\to C$ a   morphism.
Assume that  general fibers $F_{gen}$ are
\begin{enumerate}
\item smooth,
\item geometrically connected, and
\item  Fano (that is, $-K_{F_{gen}}$ is ample).
\end{enumerate}
Then every fiber $g^{-1}(c)$  contains a $k(c)$-subvariety
which is geometrically irreducible.
In particular, if every  $k(c)$-irreducible component of $g^{-1}(c)$
 is smooth (or normal), then $g^{-1}(c)$ contains a 
 $k(c)$-irreducible component 
which  is geometrically irreducible.
\end{thm}

Let us start by proving that 
 (\ref{fanofiber.thm}) $\Rightarrow$ (\ref{weakc1.thm}). 
 
Take homogeneous polynomials 
$g_1,\dots,g_s\in k[x_0,\dots,x_n]$ such that
$\deg g_i=\deg f_i$ and
$(g_1=\cdots =g_s=0)\subset \p^n_k$ is a smooth complete intersection
of dimension $n-s$.

Let $Z_1\subset \p^n_k\times \p^1_k$ be defined by the equations
$$
(uf_1+vg_1=\cdots =uf_s+vg_s=0)\subset \p^n_k\times \p^1_k,
$$
where $(u:v)$ are the coordinates on the projective line $\p^1$.
Let $Z\subset Z_1$ be the unique irreducible component
which dominates $C:=\p^1$ with projection $f:Z\to C$.

Observe that the fiber over $(0:1)$ is 
the smooth complete intersection variety 
$(g_1=\cdots =g_s=0)$, thus all but finitely many fibers are
smooth complete intersection varieties. The 
  canonical class of these fibers is the restriction of
$\o_{\p^n}(\sum \deg f_i -n-1)$, which is negative by
assumption.
Thus by (\ref{fanofiber.thm}), the fiber of $f:Z\to C$ over $(1:0)$
contains a geometrically irreducible $k$-subvariety $Y$
which is also a $k$-subvariety of $X$.\qed
\medskip

The proof of (\ref{fanofiber.thm}) proceeds in two steps.
First we use resolution of singularities to
reduce to the case $f:Y\to C$ where every fiber is a normal
crossing divisor. By this I mean that for every $c\in C$,  
every  $k(c)$-irreducible component of $f^{-1}(c)$
 is smooth. We can also achieve that 
there is an ample divisor of the form
$-mK_Y+F$ where $m>0$ and $F$ is contained in the union
of some fibers of $f$.

Then we apply  a variant of the
Koll\'ar--Shokurov Connectedness Theorem
\cite[17.4]{k-conn} to a carefully chosen
auxiliary $\q$-divisor $D$ to prove that every fiber
contains a geometrically irreducible component. One of the assumptions of
\cite[17.4]{k-conn} is, however, not satisfied in our case,
but this is compensated by other special features of the current
situation. Thus I go through the whole proof
in (\ref{conn.thm}) following 
 a quick introduction to
$\q$-divisors.

\begin{defn} Let $Y$ be a smooth $k$-variety.
A prime divisor $P$ is an irreducible and reduced codimension
1 subvariety. A $\q$-divisor is a formal linear combination
$$
D=\sum a_i P_i \qtq{where $a_i\in \q$,}
$$
and the $P_i$ are prime divisors, usually assumed
distinct.

$D$ is called a  simple normal crossing divisor
if the $P_i$ are smooth and intersect each other transversally.

The support of $D$ is $\supp D:=\cup_{i:a_i\neq 0}P_i$.
$D$ is called effective if $a_i\geq 0$ for every $i$. Every
$\q$-divisor can be uniquely written as
$$
D=D^+-D^-\qtq{where} 
D^+:=\sum_{i:a_i>0}a_iP_i\qtq{and}D^-:=\sum_{i:a_i<0}(-a_i)P_i
$$
are the positive (resp. negative) parts of $D$.
$D^+$ and $D^-$ are effective and have no common prime divisors.
Given $D=\sum a_iP_i$,  set
$$
D_{\geq 1}:=\sum_{i:a_i\geq 1}a_iP_i.
$$
A $\q$-divisor $D$ is called ample if
$mD$ is an ample (integral) divisor for some (or all)  $m>0$ such that
$ma_i$ is an integer for every $i$.

Two $\q$-divisors $D_1,D_2$ are called
$\q$-linearly equivalent (denoted by $D_1\sim_{\q} D_2$)
if there is an integer $m>0$ such that
$mD_1$ and $mD_2$ are linearly equivalent integral divisors.
(Note that even if the $D_i$ are integral divisors,
$\q$-linear equivalence is slightly different from
linear equivalence if $\pic(Y)$ contains torsion classes.)
\end{defn}

We use is the following
generalization of the Kodaira vanishing theorem.
(See  \cite[Sec.2.5]{km-book}
for a relatively short  proof and for further references.
Note that there is a misprint in the relevant
Corollary 2.68. In the last line
$\omega_Y\otimes M$ should be $\omega_Y\otimes L$.)

\begin{thm}[Kawamata--Viehweg vanishing]\label{km-vanish}
 Let $Y$ be a smooth, projective variety, $W$ any variety and
$f:Y\to W$  a  morphism.
Let  $M$ and $\Delta=\sum a_iP_i$ be $\q$-divisors
and $L$  an integral divisor 
with the following properties:
\begin{enumerate}
\item $M$ is ample,
\item $\Delta$ is a simple normal crossing divisor
and $0\leq a_i <1$ for every $i$,
\item $L\sim_{\q} M+\Delta$.
\end{enumerate}
Then $R^if_*(\o_Y(K_Y+L))=0$ for $i\geq 1$.
\end{thm}

I want to stress that 
 even in the special case when the general fiber
of $f$ is a smooth hypersurface (which is all one needs for
(\ref{main.thm.pac})) the flexibility provided by
$\q$-divisors is crucial.

Now we can prove the key technical result of this paper.

\begin{thm}[Connectedness theorem]\label{conn.thm}
 Let $Y$ be a smooth, projective  variety, $C$ a smooth curve and
$f:Y\to C$  a morphism with geometrically connected fibers.
Let $D=\sum a_i P_i$ be a 
(not necessarily effective) $\q$-divisor on $Y$ such that
\begin{enumerate}
\item $\supp D$ has simple normal crossings,
\item $D$ is $f$-vertical
(that is, its support is contained in finitely many fibers of $f$), and
\item $-(K_Y+D)$ is ample.
\end{enumerate}
Then every fiber of 
$f:\supp D_{\geq 1}\to C$ is geometrically connected.
\end{thm}

Proof. The conclusion is geometric, thus we may assume that
we are over an algebraically closed field.

For $D=D^+-D^-$,
 write $D^+=A+\Delta^+$ and $D^-=B-\Delta^-$
where $A,B$ are effective and integral, $\Delta^+$ and $\Delta^-$
are effective and if we write
$\Delta^++\Delta^-=\sum a_i P_i$ then $0\leq a_i<1$ for every $i$.
Note that $\supp A=\supp D_{\geq 1}\subset \supp D^+$
and $\supp B=\supp D^-$, thus $B$ is $f$-vertical and 
$A$ and $B$ have  no
irreducible components in common.

Consider the exact sequence
$$
0\to \o_Y(B-A)\to \o_Y(B)\to \o_A(B)\to 0,
$$
and apply $f_*$  to get
$$
f_*\o_Y(B)\to f_*\o_A(B)\to R^1f_*\o_Y(B-A).
$$
Since $B$ is vertical and $f$ has geometrically connected fibers,
$f_*\o_Y(B)$ is a torsion free sheaf which is generically
isomorphic to $\o_C$. Thus $f_*\o_Y(B)$ is a rank 1 locally free sheaf
on $C$.
Observe that
$$
\begin{array}{rcl}
B-A&\sim_{\q}& D^-+\Delta^--D^++\Delta^+\sim_{\q} -D+( \Delta^-+\Delta^+)\\
&\sim_{\q}& K_Y+\bigl(-(K_Y+D)\bigr)+( \Delta^-+\Delta^+),
\end{array}
$$
hence by (\ref{km-vanish}) we conclude that $R^1f_*\o_Y(B-A)=0$ and
so $f_*\o_A(B)$ is the quotient of the 
rank 1 locally free sheaf $f_*\o_Y(B)$.

$B$ is an effective divisor which has no
irreducible components in common with $A$, thus
we have an injection $\o_A\into \o_A(B)$
and so an injection
$f_*\o_A\into f_*\o_A(B)$.

Therefore $f_*\o_A$  is the quotient of a subsheaf $M$ of  the 
rank 1 locally free sheaf $f_*\o_Y(B)$.
Since $C$ is a smooth curve, its local rings $\o_{c,C}$ are
principal ideal domains, hence $M$ itself is a
rank 1 locally free sheaf.

Finally, set  $A(c):=A\cap f^{-1}(c)$ and let $A(c)_1,\dots ,A(c)_m$
be its connected components. 
We have surjections
$$
\o_{c,C}\cong M_c\onto f_*\o_{A(c)}\cong H^0(A(c), \o_{A(c)})
\cong \sum_{i=1}^m H^0(A(c)_i, \o_{A(c)_i})\onto k^m.
$$
This implies that $m=1$, hence
 $f|_A$ has 
geometrically connected fibers.\qed

\begin{cor}\label{geomirredcomp.cor}
  Let $Y$ be a smooth, projective  variety, $C$ a smooth curve and
$f:Y\to C$  a morphism with geometrically connected fibers.
Let $D=\sum a_i P_i$ be a (not necessarily effective)  
$\q$-divisor on $Y$ such that
\begin{enumerate}
\item $\supp D+(\mbox{any fiber of $f$})$ has simple normal crossings,
\item $D$ is $f$-vertical,
\item $-(K_Y+D)$ is ample.
\end{enumerate}
Them, for every $c\in C$, the fiber $F_c$ contains
a $k(c)$-irreducible component which is
geometrically irreducible.
\end{cor}

Proof.  Given $c\in C$, let $A$ be any  divisor on $C$ which
is linearly equivalent to 0 such that $c\in \supp A$.
Let $F_A:=f^{-1}(A)$. 
$F_A$ is linearly equivalent to 0, thus we can
add any rational multiple of $F_A$ to $D$ without changing the
assumptions of (\ref{conn.thm}).
Hence we may assume that
\begin{enumerate}\setcounter{enumi}{3}
\item in a neighborhood of $F_c$, every
irreducible component of $D$ has coefficient $\leq 1$, and
\item at least one irreducible component of $F_c$ has coefficient $1$ in $D$.
\end{enumerate}
Let $E\subset F_c$ be such a component. We claim that
$E$ is geometrically irreducible.

To see this, let $m$ be the multiplicity of $E$ in $F_A$ and
consider $D'=\sum a'_iP_i$ defined by 
$$
D':=D-\tfrac{\epsilon}{m}F_A+\epsilon E\qtq{for $0<\epsilon\ll 1$.}
$$
This choice assures that  
in a neighborhood of $F_c$, $E$ is the only 
irreducible component of $F_c$ which has coefficient $1$ in $D'$.
Furthermore
$$
-(K_Y+D')\sim_{\q} -(K_Y+D)-\epsilon E
$$
is still ample for $0<\epsilon\ll 1$
since ampleness is an open condition by Kleiman's criterion
\cite{klei}.
Thus, in a neighborhood of $F_c$, $E=D'_{\geq 1}$ and so
$E$ is  a smooth and geometrically connected variety, hence
geometrically irreducible. \qed

\begin{rem} The assumption  (\ref{geomirredcomp.cor}.2)
can be relaxed considerably. The proof works without changes
if (\ref{geomirredcomp.cor}.2)  is replaced by the following two conditions:
\begin{enumerate}
\item[(2.i)] $D_{\geq 1}$ is $f$-vertical, and
 \item[(2.ii)] $H^0(F_{gen}, \o_{F_{gen}}(\rup{D^-}))$ is 1-dimensional
over $k(C)$, where $F_{gen}$ is the generic fiber of $f$.
\end{enumerate}
The general case is also useful since it can be used to prove that
(\ref{fanofiber.thm}) also holds when the generic fiber is a
$\q$-Fano variety
\end{rem}
\medskip

\medskip

Proof of (\ref{fanofiber.thm}). Choose an embedding of the generic
fiber $F_{gen}\into \p^N$ such that
$-mK_{F_{gen}}\sim H|_{F_{gen}}$ for some $m\geq 1$
where  $H$ is the hyperplane class on $\p^N$.
Let $B$ be a very ample divisor on $C$. Then 
$\pi_C^*B+\pi_P^*H$ is very ample on $\p^N\times C$
where $\pi_C,\pi_P$ are the coordinate projections.

By taking the closure of $F_{gen}$ we get a variety
  $V\subset \p^N\times C$ with a projection map
$\pi:V\to C$  and a birational map $\phi:V\map Z$
which is an isomorphism between   
 $V^0:=\pi^{-1}(C^0)$ and $Z^0:=f^{-1}(C^0)$
for some open $C^0\subset C$.
Moreover, by shrinking $C^0$ if necessary, we may also
assume that $V^0\to C^0$ is smooth and  
$$
-mK_V|_{V^0}\sim (\pi_C^*B+\pi_P^*H)|_{V^0}.
$$
We can apply Hironaka's resolution theorem
to get $h_1:V_1\to V$ such that $V_1$ is smooth and  
$\phi$ lifts to a morphism
$\phi_1:V_1\to Z$.

Let $F^{sing}\subset V_1$ be the union of all singular fibers
of $\pi\circ h_1$. 
Applying Hironaka's resolution theorem again to
$F^{sing}\subset V_1$ we get
$h_2:V_2\to V_1$ and $h:V_2\to V$
such that every fiber of $\pi\circ h:V_2\to C$
is a simple normal crossing divisor.
That is, for every $c\in C$.
every $k(c)$-irreducible component of
the reduced fiber $\red F_c$ is smooth.

Throughout these resolutions we do not blow up anything
above $C^0$.

Since $h:V_2\to V$ is a composite of blow ups
of subvarieties, there is an $m_2>0$ and an
$h$-exceptional divisor $E$ such that
$h^*(m_2(\pi_C^*B+\pi_P^*H))-E$ is ample on $V_2$
(cf.\ \cite[II.7.10.b and II.7.13]{harts}). Dividing by $m\cdot m_2$
we conclude that there is an ample  $\q$-divisor $M$ on $V_2$
such that
$$
-K_{V_2}|_{V^0}\sim_{\q} M|_{V^0}.
$$
Thus there is a $\q$-divisor $D$ supported
in $V_2\setminus V^0$ such that
$$
-(K_{V_2}+D)\sim_{\q} M.
$$
Since the support of $D$ is contained in a union
of fibers of $\pi\circ h$, it is a 
simple normal crossing divisor. Thus
$D$ is vertical and the assumptions
(\ref{geomirredcomp.cor}.1--3) hold.

Hence by  (\ref{geomirredcomp.cor}), every fiber of $\pi\circ h:V_2\to C$
contains a geometrically irreducible component.

Since every fiber of $f:Z\to C$ is dominated by a fiber 
of $\pi\circ h:V_2\to C$, we conclude that
every fiber of $f:Z\to C$
contains a geometrically irreducible subvariety.

Finally, assume that  every  $k(c)$-irreducible component of $g^{-1}(c)$
 is smooth (or normal).
Let $W\subset g^{-1}(c)$ be a geometrically irreducible subvariety
and $F\subset g^{-1}(c)$ an irreducible
component containing $W$. Write $F_{\bar k}=F_1+\cdots+F_m$
where the $F_i$ are irreducible over $\bar k$.
One of the $F_i$ contains $W_{\bar k}$, but then so do all the others since
the $F_i$ are conjugate over $k$. Since
$F_{\bar k}$ is normal, this implies that $m=1$ and
$F$ is geometrically irreducible.
\qed

\medskip

Theorem (\ref{fanofiber.thm}) naturally raises the
following question:

\begin{ques} Which ``natural'' classes of
schemes $\mathbb S$ satisfy the following
property
$$
\begin{array}{c}
\mbox{For every field $k$ and for every $k$-scheme
$S\in \mathbb S$,}\\
\mbox{$S$ contains a geometrically irreducible subscheme.}
\end{array}
\eqno{(*)}
$$
\end{ques}

We have shown that $(*)$  holds for
$$
{\mathbb S}=\{\mbox{degenerations of Fano varieties}\}.
$$

There are two immediate generalizations, but
$(*)$ fails for both.
First, degenerations of Fano varieties are all
rationally chain connected, that is, any two points are
connected by a chain of rational curves.
(See \cite[Chap.IV]{rcbook} for a general overview.)

The triangle $(xyz=0)\subset \p^2$ is
rationally chain connected. Let $K/\q$ be any cubic extension
with norm form $N(x,y,z)$. Then $C_N:=(N(x,y,z)=0)$ has no
geometrically irreducible $\q$-subvarieties
but it is isomorphic to the triangle over $\bar \q$.

One can also try to work with singular Fano schemes.
That is, schemes $X$ such that $\omega_X$ is a line bundle
such that $\omega_X^{-1}$ is ample.
Here $(*)$ again fails.

Take the affine variety $(N(x,y,z)+x^4+y^4+z^4=0)\subset \a^3$.
Blow up the origin to get $Y$. The exceptional curve is
isomorphic to $C_N$, let $I$ be its ideal sheaf.
Then $X=\spec_Y \o_Y/I^2$ is a Fano scheme
 with no
geometrically irreducible $\q$-subvarieties.

I have, however, no counter examples to the following questions:

\begin{ques} Does $(*)$ hold for the following
two classes of schemes:
\begin{enumerate}
\item Degenerations of smooth rationally  connected varieties.
\item Reduced Fano schemes.
\end{enumerate}
\end{ques}

In fact, in both cases it may be true that such a scheme
 contains a geometrically irreducible subscheme
which is also rationally connected.
The following variant is especially interesting:

\begin{ques}\label{src.questions}
 Let $k$ be a field of characteristic 0, $C$ a smooth $k$-curve,
$Z$ a smooth $k$-variety and
$g:Z\to C$ a  projective  morphism.
Assume that
\begin{enumerate}
\item  the generic fiber $F_{gen}$ is  rationally connected,
\item every fiber is a normal crossing divisor, and
\item for every $c\in C$, 
every  $k(c)$-irreducible component of $g^{-1}(c)$
 is smooth.
\end{enumerate}
Is it true that  every  fiber
$g^{-1}(c)$ contains a 
 $k(c)$-irreducible component which  
 is  rationally connected (and hence geometrically irreducible)?
\end{ques}

\begin{rem}[Positive characteristic]
The conjecture of Ax needs only minor
modifications in positive characteristic, see
\cite[Chap.21]{f-j}.

It is known that for any prime $p$, the following are equivalent:
 \begin{enumerate}
\item $\f_p(t)$ is weakly  $C_1$.
\item Every  field of characteristic p is  weakly  $C_1$.
\item Every perfect PAC field of characteristic $p$ is $C_1$.
\end{enumerate}

My proof has 2 difficulties in
general. First,  resolution of singularities is not known
 (but it is expected to be true).
Second, Kodaira's vanishing theorem and its generalization
(\ref{km-vanish})
are  false in positive characteristic. As far as I know,
however, the Koll\'ar--Shokurov Connectedness Theorem
may hold in positive characteristic.
\end{rem}

\begin{ack}  I thank M.\ Jarden, J.\ McKernan and T.\ Szamuely
for useful comments and suggestions.
Partial financial support was provided by  the NSF under grant number 
DMS-0500198. 
\end{ack}

\bibliography{refs}

\providecommand{\bysame}{\leavevmode\hbox to3em{\hrulefill}\thinspace}
\providecommand{\MR}{\relax\ifhmode\unskip\space\fi MR }
\providecommand{\MRhref}[2]{%
  \href{http://www.ams.org/mathscinet-getitem?mr=#1}{#2}
}
\providecommand{\href}[2]{#2}
\begin{thebibliography}{DJL83}

\bibitem[Ax68]{ax}
James Ax, \emph{The elementary theory of finite fields}, Ann. of Math. (2)
  \textbf{88} (1968), 239--271. \MR{MR0229613 (37 \#5187)}

\bibitem[DJL83]{djl}
Jan Denef, Moshe Jarden, and D.~J. Lewis, \emph{On {A}x-fields which are
  {$C\sb{i}$}}, Quart. J. Math. Oxford Ser. (2) \textbf{34} (1983), no.~133,
  21--36. \MR{MR688420 (84j:12028)}

\bibitem[FJ05]{f-j}
Michael~D. Fried and Moshe Jarden, \emph{Field arithmetic}, Ergebnisse der
  Mathematik und ihrer Grenzgebiete. 3. Folge. A Series of Modern Surveys in
  Mathematics [Results in Mathematics and Related Areas. 3rd Series. A Series
  of Modern Surveys in Mathematics], vol.~11, Springer-Verlag, Berlin, 2005.
  \MR{MR2102046 (2005k:12003)}

\bibitem[Har77]{harts}
Robin Hartshorne, \emph{Algebraic geometry}, Springer-Verlag, New York, 1977,
  Graduate Texts in Mathematics, No. 52. \MR{57 \#3116}

\bibitem[Kle66]{klei}
Steven~L. Kleiman, \emph{Toward a numerical theory of ampleness}, Ann. of Math.
  (2) \textbf{84} (1966), 293--344. \MR{MR0206009 (34 \#5834)}

\bibitem[KM98]{km-book}
J{\'a}nos Koll{\'a}r and Shigefumi Mori, \emph{Birational geometry of algebraic
  varieties}, Cambridge Tracts in Mathematics, vol. 134, Cambridge University
  Press, Cambridge, 1998. \MR{MR1658959 (2000b:14018)}

\bibitem[Kol92]{k-conn}
J{\'a}nos Koll{\'a}r, \emph{Adjunction and discrepancies}, Flips and abundance
  for algebraic threefolds, Soci\'et\'e Math\'ematique de France, Paris, 1992,
  pp.~183--192. \MR{MR1225842 (94f:14013)}

\bibitem[Kol96]{rcbook}
\bysame, \emph{Rational curves on algebraic varieties}, Springer-Verlag,
  Berlin, 1996. \MR{98c:14001}

\end{thebibliography}

\bigskip

\noindent Princeton University, Princeton NJ 08544-1000

\begin{verbatim}kollar@math.princeton.edu\end{verbatim}

\end{document}